\documentclass[11pt,letterpaper]{amsart}
\usepackage{amsmath,amssymb}
\usepackage{mathrsfs}
\usepackage{amsfonts}
\usepackage[dvipsnames]{xcolor}
\usepackage{graphicx}
\usepackage{float}  
\usepackage{cite}
\usepackage{microtype}
\usepackage{cases}
\usepackage{comment}
\usepackage{indentfirst}
\allowdisplaybreaks
\theoremstyle{plain}
\newtheorem{theorem}{Theorem}[section]
\newtheorem{lemma}[theorem]{Lemma}
\newtheorem{proposition}[theorem]{Proposition}

\theoremstyle{definition}

\theoremstyle{remark}

\numberwithin{equation}{section}

\DeclareMathOperator{\tr}{tr}
\DeclareMathOperator{\dist}{dist}
\newcommand{\R}{\mathbb R}
\newcommand{\Om}{\Omega}
\newcommand{\eps}{\varepsilon}
\newcommand{\Ga}{\Gamma_k}
\newcommand{\cD}{\mathcal D}
\newcommand{\cK}{\mathcal K}
\newcommand{\cL}{\mathcal L}
\newcommand{\dm}{\,d\mu}

\usepackage[colorlinks=true,citecolor=blue,linkcolor=blue]{hyperref}
\hypersetup{pdftitle={Global $W^{2,p}$ and $C^{2,alpha}$ regularity for the $k$-Hessian equation},
  pdfauthor={Yasheng Lyu}}

\begin{document}
\title[Global $W^{2,p}$ and $C^{2,\alpha}$ regularity]
{Global $W^{2,p}$ and $C^{2,\alpha}$ regularity for the $k$-Hessian equation}

\author[Y. Lyu]{Yasheng Lyu}
\address{School of Mathematics and Statistics, Xi'an Jiaotong University, Xi'an, Shaanxi 710049, People's Republic of China}
\email{lvysh21@stu.xjtu.edu.cn}

\date{}
\keywords{$k$-Hessian equation, global $W^{2,p}$ estimates, global $C^{2,\alpha}$ estimates, Hessian measures, viscosity solutions.}

\begin{abstract}
In this paper, we prove that for $3\leq k\leq n-1$, continuous $k$-convex
solutions of the $\sigma_{k}$ equation on bounded strictly $(k-1)$-convex
$C^{3}$ domains $\Omega\subset\mathbb{R}^{n}$ with $C^{3}$ boundary values and positive
continuous density on $\overline{\Omega}$ belong to $W^{2,p}(\Omega)$ for
every $1<p<\infty$.
If, in addition, the density is $C^{\alpha}$ on $\overline{\Omega}$,
$0<\alpha<1$, the solutions belong to $C^{2,\alpha}(\overline{\Omega})$.
\end{abstract}

\maketitle

\section{Introduction}\label{sec:intro}

We consider the Dirichlet problem
\begin{equation}\label{eq:dirichlet}
 \sigma_k(D^2u)=f\quad\text{in }\Om,\qquad
 u=\varphi\quad\text{on }\partial\Om,
\end{equation}
where $\Om\subset\R^n$ is a bounded domain and $3\le k\le n-1$.

For $\lambda\in\R^m$ and $1\le j\le m$, let
\[
 \sigma_j(\lambda)=\sum_{1\le i_1<\cdots<i_j\le m}
  \lambda_{i_1}\cdots\lambda_{i_j},\qquad
 \Gamma_j=\{\lambda\in\R^m:\sigma_i(\lambda)>0,\ 1\le i\le j\},
\]
with $\sigma_0=1$. For a symmetric matrix $A$ with eigenvalue vector
$\lambda(A)$, we write $\sigma_j(A)=\sigma_j(\lambda(A))$ and
$A\in\Gamma_j$ if $\lambda(A)\in\Gamma_j$.
A $C^2$ function is $k$-admissible if $D^2u\in\Ga$.

Let $\kappa=(\kappa_1,\ldots,\kappa_{n-1})$ denote the principal
curvatures of $\partial\Om$ relative to the interior unit normal $\nu$,
with the sign convention that spheres have positive curvatures. The
boundary is strictly $j$-convex if $\kappa(x)\in\Gamma_j$ at every $x\in\partial\Om$.
For $j=k-1$ this is equivalent to condition~(10) of Caffarelli,
Nirenberg and Spruck~\cite{CNS}, which requires a fixed $R>0$ such that
$(\kappa(x),R)\in\Ga$ at every $x\in\partial\Om$.

We denote by $\Phi_k(\Om)$ the class of $k$-convex functions. A
continuous function is $k$-convex if every $C^2$ function touching it
from above has Hessian in $\overline\Ga$ at the point of contact.
We interpret the equation in \eqref{eq:dirichlet} as
$\mu_k[u]=f\,dx$, where $\mu_k[u]$ is the Hessian measure of Trudinger
and Wang~\cite{TW1,TW2}. For positive continuous $f$ this agrees with
the viscosity formulation~\cite[Section~2]{TW2}.

\begin{samepage}
\begin{theorem}\label{thm:wp}
Let $3\le k\le n-1$ and let $\Om\subset\R^n$ be a bounded strictly
$(k-1)$-convex domain with $\partial\Om\in C^3$.
Suppose $\varphi\in C^3(\partial\Om)$ and $f\in C(\overline\Om)$,
with $0<\lambda\le f\le\Lambda$.
Then any $u\in C(\overline\Om)\cap\Phi_k(\Om)$ solving
\eqref{eq:dirichlet} in the Hessian-measure sense belongs to $W^{2,p}(\Om)$ for every $1<p<\infty$
and satisfies
\begin{equation}\label{eq:main-estimate}
 \|u\|_{W^{2,p}(\Om)}\le C_p,
\end{equation}
where $C_p$ depends only on $n,k,p,\lambda,\Lambda,\Om,\varphi$
and the modulus of continuity of $f$.
\end{theorem}
\end{samepage}

\begin{theorem}\label{thm:c2alpha}
Let $3\le k\le n-1$ and let $\Om\subset\R^n$ be a bounded strictly
$(k-1)$-convex domain with $\partial\Om\in C^3$.
Suppose $\varphi\in C^3(\partial\Om)$ and $f\in C^\alpha(\overline\Om)$
for some $0<\alpha<1$, with $f\ge\lambda>0$.
Then any $u\in C(\overline\Om)\cap\Phi_k(\Om)$ solving
\eqref{eq:dirichlet} in the Hessian-measure sense belongs to $C^{2,\alpha}(\overline\Om)$ and satisfies
\begin{equation}\label{c2:eq:main-estimate}
 \|u\|_{C^{2,\alpha}(\overline\Om)}\le C,
\end{equation}
where $C$ depends only on $n,k,\alpha,\lambda$,
$\|f\|_{C^\alpha(\overline\Om)},\Om$ and $\varphi$.
\end{theorem}

The dependence on $\Om$ and $\varphi$ in both estimates is through the
diameter of $\Om$, the $C^3$ norms of the boundary charts and boundary
values, the moduli of continuity of their third derivatives, and lower
bounds for the chart radii and for $\dist(\kappa,\partial\Gamma_{k-1})$
on $\partial\Om$.

For smooth data with $f>0$, Caffarelli, Nirenberg and
Spruck~\cite[Theorem~3]{CNS} proved existence and uniqueness of a
smooth admissible solution when $\partial\Om$ is connected and
$\sigma_{k-1}(\kappa)>0$ on $\partial\Om$.
Ivochkina~\cite{Ivochkina} treated $k=3$ with zero boundary values
on strictly convex domains. Guan~\cite{Guan} proved solvability in the
presence of a strict admissible subsolution. On strictly $(k-1)$-convex
$C^{3,1}$ domains with $C^{3,1}$ boundary values,
Trudinger~\cite[Theorem~1.1]{Trudinger1995} obtained classical
solvability for positive $C^{1,1}$ density. The global Hessian estimate of Ivochkina,
Trudinger and Wang~\cite{ITW} involves second derivatives of $f^{1/k}$
and is independent of $\inf f$.

For the Monge--Amp\`ere equation, Caffarelli~\cite{Caffarelli} proved
interior $W^{2,p}$ estimates for strictly convex solutions with positive
continuous density. Savin~\cite[Theorem~1.2 and Corollary~1.4]{Savin}
proved the global estimate under quadratic separation at the boundary,
in particular on uniformly convex $C^3$ domains with $C^3$ boundary
values. Under the latter boundary assumptions, Trudinger and
Wang~\cite[Theorem~1.1]{TWMA} obtained global $C^{2,\alpha}$ regularity
for positive $C^\alpha$ density. Savin~\cite[Theorem~1.1]{SavinPoint}
proved a pointwise boundary $C^{2,\alpha}$ estimate.

For $k=2$, the author proved interior $W^{2,p}$ regularity for
positive continuous density~\cite{LyuInterior} and the counterparts of
Theorems~\ref{thm:wp} and~\ref{thm:c2alpha} on strictly mean-convex
$C^3$ domains with $C^3$ boundary values~\cite{LyuK2}. For positive $C^\alpha$ density, Zhou and
Zhu~\cite{ZZ} proved interior $C^{2,\alpha}$ regularity for convex
viscosity solutions, and Chen, Zhou and Zhu~\cite{CZZ} proved it for
admissible solutions. For densities bounded above and below by positive
constants, the author~\cite{LyuEps} proved that smooth solutions on
strictly mean-convex $C^{2,\beta}$ domains with $C^{2,\beta}$ boundary
values satisfy a global $W^{2,1+\eps}$ estimate for some $\eps>0$. For boundary
counterexamples with $C^{2,1}$ boundary or boundary values, see
Wang~\cite{Wang} for $n=2$ and~\cite[Section~4]{LyuK2} for $k=2$.

For each $3\le k\le n-1$, the extension of Mooney's
example~\cite[Sections~5 and~6]{MooneySingular}
in~\cite[Proposition~3.1]{LyuWone} gives a convex solution $u$ of
$\sigma_k(D^2u)=1$ in a ball such that $u\notin W^{2,p}_{\mathrm{loc}}$
for any $p>1$. For earlier singular examples see Urbas~\cite{Urbas}.
For positive $C^{1,1}$ density and zero
boundary values, Chou and Wang~\cite[Theorem~1.5]{CW} proved
$|u|^4|D^2u|\le C$, with $C$ depending also on $\|u\|_{C^1}$.
For $0<\lambda\le f\le\Lambda$, Chen, Wang and Wu~\cite{CWW} proved
interior $C^{1,\beta}$ regularity for convex solutions with zero
boundary values on bounded convex domains. For smooth positive density, Urbas obtained
interior Hessian bounds from $W^{2,p}$ bounds, first for
$p>kn/2$~\cite{Urbas2000} and then for $p>k(n-1)/2$~\cite{Urbas2001}.
For $k<n/2$, Dinew and Ko{\l}odziej~\cite[Remark~4.4]{DK} observed
that $p>n(k-1)/2$ suffices.

The rest of the paper is organized as follows.
Section~\ref{sec:wp} proves Theorem~\ref{thm:wp}.
Section~\ref{sec:c2alpha} proves Theorem~\ref{thm:c2alpha}.

\section{Global \texorpdfstring{$W^{2,p}$}{W2,p} estimates}\label{sec:wp}

Let $G=\sigma_k^{1/k}$ and fix $0<m<M$. We first consider smooth data.
The constants $c,C$ may depend on $n,k,m,M$ and on the domain and
boundary values in Theorem~\ref{thm:wp}. We write $C_b$ if the constant
also depends on the $C^2$ norm of a smooth function $b$ approximating
$f^{1/k}$ from below.

For a function $F$ of a symmetric matrix, we use the symmetric
coefficients $F^{ij}$ and $F^{ij,rs}$ determined by
\[
 \begin{aligned}
 \left.\frac{d}{dt}F(A+tH)\right|_{t=0}
    &=F^{ij}(A)H_{ij},\\
 \left.\frac{\partial^2}{\partial t\,\partial s}F(A+tH+sE)\right|_{t=s=0}
    &=F^{ij,rs}(A)H_{ij}E_{rs}
 \end{aligned}
\]
for symmetric $H,E$. We use the Frobenius norm $|A|$.
Repeated indices are summed over all their
values, including both orders of each off-diagonal pair. In an
orthonormal frame adapted to a hyperplane, Greek indices range over
the $n-1$ tangential directions.

\subsection{The tangential Hessian at the boundary}
\label{wp:nondegeneracy}

Fix a $C^3$ extension $\Phi$ of $\varphi$ to a neighborhood of
$\overline\Om$. In each boundary chart extend the boundary value
constantly in the vertical direction and combine these extensions by a
smooth partition of unity. This gives bounds for $\|\Phi\|_{C^3}$ and
the modulus of $D^3\Phi$ in terms of the boundary data, and gives a
smooth extension for smooth data. Write $r(x)=\dist(x,\partial\Om)$ for
$x\in\overline\Om$.

\begin{lemma}\label{bdy:barriers}
There exist a $C^3$ defining function $\rho$ and $\theta>0$ such that,
with $d=-\rho$,
\begin{equation}\label{eq:rho}
d\asymp r,\qquad D^2\rho-\theta I\in\Ga
\quad\text{on }\overline\Om.
\end{equation}
Near $\partial\Om$ one has $\rho=-r+Jr^2$ for some $J>0$. If $v$ is smooth, $k$-admissible,
$v=\varphi$ on $\partial\Om$, and $m\le G(D^2v)\le M$, then
\[
 |v-\Phi|\le C d,\qquad
 \|v\|_{L^\infty(\Om)}+\|Dv\|_{L^\infty(\overline\Om)}\le C
\]
where $C$ is uniform for the above data.
\end{lemma}

\begin{proof}
Choose principal coordinates at a boundary point, with $e_n=\nu$.
In this frame the derivatives of $r$ from the interior satisfy
\cite[Section~2]{CNS}
\[
 Dr=\nu,\qquad D^2r=\operatorname{diag}(-\kappa_1,\ldots,-\kappa_{n-1},0).
\]
Thus $D^2(-r+Jr^2)$ has eigenvalues
$(\kappa_1,\ldots,\kappa_{n-1},2J)$. Since
$\sigma_j(\kappa,2J)=\sigma_j(\kappa)+2J\sigma_{j-1}(\kappa)$
and $\kappa\in\Gamma_{k-1}$ uniformly, this Hessian lies in $\Ga$
for large $J$. Fix such a $J$ and a small $\ell>0$ so that the inclusion
holds on $0\le r\le \ell$ and $J\ell<1/8$. In this strip take a smooth
regularized maximum of $-r+Jr^2$ and a strictly convex quadratic with
values in $(-\ell/2,-\ell/4)$ on $\overline\Om$. For sufficiently small
regularization parameter, the maximum agrees with $-r+Jr^2$ near
$r=0$ and with the quadratic near $r=\ell$. Extend it by the quadratic
to the rest of $\Om$. Its Hessian belongs to $\Ga$ by convexity and
ellipticity of the cone. Compactness now gives \eqref{eq:rho}. By
concavity and homogeneity of $G$, for $a^{ij}=G^{ij}(A)$ with $A\in\Ga$,
\begin{equation}\label{eq:linear-barrier}
 a^{ij}d_{ij}\le-\theta\tr a,\qquad \tr a\ge G(I).
\end{equation}

For large $L$, the function $\Phi+L\rho$ is admissible and
$G(D^2(\Phi+L\rho))\ge M$, while $\Delta(\Phi-L\rho)\le0$.
Comparison gives $\Phi+L\rho\le v\le\Phi-L\rho$. Given $x\in\Om$,
take $R=r(x)/2$ and a nearest boundary point $y$. The barriers and
the Lipschitz continuity of $\Phi$ give
$\sup_{B_R(x)}|v-\varphi(y)|\le CR$. Since
$\|\sigma_k(D^2v)\|_{L^\infty(B_R(x))}^{1/k}\le M$, the gradient estimate
of Liu, Ma and Zang~\cite[Corollary~1.4]{LMZ}, with $q=\infty$
and affine function equal to $\varphi(y)$, gives
$|Dv(x)|\le C(R^{-1}\sup_{B_R(x)}|v-\varphi(y)|+RM)\le C$.
For the corresponding gradient potential estimate, see also
He~\cite[Theorem~1.1]{He}.
\end{proof}

\begin{proposition}\label{bdy:margin}
Suppose $v\in C^2(\overline\Om)$ is $k$-admissible,
$v=\varphi$ on $\partial\Om$, and $G(D^2v)\ge m>0$ in $\Om$.
There exists $c_0>0$, depending only on $n,k,m,\Om$ and $\varphi$, such that
\begin{equation}\label{eq:tangential-margin}
 D^2v|_{\nu^\perp}-c_0 I\in\Gamma_{k-1}
 \qquad\text{on }\partial\Om .
\end{equation}
In particular, $\sigma_{k-1}(D^2v|_{\nu^\perp})$ has a uniform
positive lower bound. If also $G(D^2v)\le M$, the norm of the
tangential Hessian has a uniform upper bound depending in addition on $M$.
\end{proposition}

\begin{proof}
We use the supporting-plane argument of Caffarelli, Nirenberg and
Spruck~\cite[Section~6]{CNS} and the $C^3$ barrier
of~\cite[Proposition~2.1]{LyuK2}. For the tangential estimate in the
Monge--Amp\`ere case, see Wang~\cite[Lemma~3.1]{Wang}.
Center orthonormal coordinates at the boundary point, with $e_n=\nu$
and the other axes in principal directions. Write $x_n=\beta(x')$ for
the boundary, with $\Om$ above the graph, $\beta(0)=0$ and $D'\beta(0)=0$.
Set $\phi(x')=\varphi(x',\beta(x'))$. Then
$K=D^2_{x'}\beta(0)=\operatorname{diag}(\kappa_1,\ldots,\kappa_{n-1})$.
Strict $(k-1)$-convexity gives $K-2\gamma I\in\Gamma_{k-1}$ for
some uniform $\gamma>0$. Let $q$ be the upper endpoint of
$\{t:D^2_{x'}\phi(0)-tK\in\Gamma_{k-1}\}$. This interval contains
$(-\infty,-C)$ by the bound on $D^2_{x'}\phi$ and the inclusion for $K$.
Its upper endpoint is uniformly bounded by the trace condition.
Subtracting $\phi(0)+D'\phi(0)\cdot x'+qx_n$ from $v$ and its
boundary values, we may thus assume that
\[
 \phi(0)=0,\qquad D'\phi(0)=0,\qquad
 D^2_{x'}\phi(0)\in\partial\Gamma_{k-1}.
\]
Choose a supporting functional $B\mapsto\tr(PB)$ at $D^2_{x'}\phi(0)$
which is nonnegative on $\Gamma_{k-1}$. Since the closed cone
contains all positive semidefinite matrices, $P\ge0$. Normalize $\tr P=1$.
Then $\tr(PD^2_{x'}\phi(0))=0$ and $\tr(PK)\ge2\gamma$.

Let $s=|x'|$, and let $\beta_2,\phi_2,\phi_3$ denote the homogeneous
terms of the indicated degrees in the Taylor expansions at the origin.
Set
\[
 d_0=x_n-\beta_2+\gamma s^2/2,\qquad
 w=\phi_2+\phi_3+(l\cdot x')(x_n-\beta_2),
\]
where $(\tr(PK)I+2\gamma P)l=D'\tr(PD^2_{x'}\phi_3)$.
In particular, $|l|\le C$. On $x_n=\beta(x')$ one has
$d_0=\gamma s^2/2+O(s^3)$ and $|w-\phi|\le\eps(h)s^3$ for $s<h$,
with $\eps(h)\to0$ uniformly in the boundary point. In
$U_h=\{s<h,\ \beta(x')<x_n<h^2\}$ one has $d_0\asymp s^2+|x_n|$
for small $h$, since $d_0\ge x_n-\beta(x')+\gamma s^2/4$.

For $i<n$, write $X_i=e_i-(d_0)_i e_n$ and set
$L\psi=\sum_{i,j<n}P_{ij}\,X_i\cdot(D^2\psi\,X_j)$. Then $Ld_0=-\tr(PK)+\gamma$
and $Lw=0$. Indeed, $X_i\cdot D(x_n-\beta_2)=-\gamma x_i$, so
\[
 Lw=\tr(PD^2_{x'}\phi_3)-(l\cdot x')\tr(PK)-2\gamma(Pl)\cdot x'=0 .
\]
Choose $C_0$ so that $C_0s^2-\beta(x')\ge s^2$ and put
\[
 V=w+Nd_0^2/2+\eta(C_0s^2-x_n),\qquad N\ge1,\quad 0<\eta\le1 .
\]

Suppose $D^2V\in\Ga$ at a point of $U_h$. Its restriction to
$(Dd_0)^\perp$ belongs to $\Gamma_{k-1}$. At this point take an
orthonormal basis of that hyperplane of the form $X_i+O(s^2)$, which
is possible since $X_i\cdot X_j=\delta_{ij}+O(s^2)$.
Let $B$ be the matrix of the restriction in this basis. The term
$N Dd_0\otimes Dd_0$ vanishes on the hyperplane. Since $D^2w$
and $D^2d_0$ are bounded, the two bases give
\[
 0<\tr(PB)=LV+O(s^2)(1+Nd_0)
 \le C\eta+Cs^2-(\gamma-Cs^2)Nd_0 .
\]
Fix $h$ small and $N$ large. As $s^2\le C d_0$, this gives
$Nd_0\le C\eta$, $|B|\le C$, and $\tr(PB)\le C\eta$.
Since $\tr P=1$, $B-tI$ reaches $\partial\Gamma_{k-1}$ for some
$0<t\le\tr(PB)$. There $\sigma_{k-1}$ vanishes. Its Lipschitz
bound on bounded sets now gives $\sigma_{k-1}(B)\le C\eta$.
If $\sigma_k(B)>0$, then $B\in\Ga$ and Maclaurin's inequality gives
$\sigma_k(B)\le C\eta^{k/(k-1)}$. Writing $e=Dd_0/|Dd_0|$ and
$z$ for the mixed column, we obtain
\[
 \sigma_k(D^2V)
 =\sigma_k(B)+V_{ee}\sigma_{k-1}(B)-\sigma_{k-1}^{\alpha\beta}(B)z_\alpha z_\beta
 \le C\bigl(N\eta+\eta^{k/(k-1)}\bigr),
\]
since $V_{ee}\le CN$ and $(\sigma_{k-1}^{\alpha\beta}(B))$ is positive definite.

Take $\eta=\delta/N$ with $\delta>0$ fixed so small that the last
bound is less than $m^k$. On $x_n=\beta(x')$,
\[
 V-\phi\ge cNs^4+\eta s^2-\eps(h)s^3\ge0
\]
after decreasing $h$ so that $\eps(h)\le2\sqrt{c\delta}$.
On the remaining boundary of $U_h$ one has $d_0\ge ch^2$.
Subharmonicity bounds $v$ from above by its boundary maximum, so
increasing $N$ gives $V\ge v$ there as well. At a positive interior
maximum of $v-V$, one would have $D^2V\ge D^2v$ and hence
$D^2V\in\Ga$ with $\sigma_k(D^2V)\ge m^k$. Thus $v\le V$.
At the origin equality holds, and $v_n(0)\le V_n(0)=-\eta$.

Twice differentiating $v(x',\beta(x'))=\phi(x')$ gives
$D^2_{x'}v(0)=D^2_{x'}\phi(0)-v_n(0)K$. Since $K-2\gamma I\in\Gamma_{k-1}$,
it follows that $D^2_{x'}v(0)-\eta\gamma I\in\Gamma_{k-1}$. The choices
of $h,N,\delta$ are uniform under the stated bounds and moduli of
continuity. When $G(D^2v)\le M$, comparison with $\Phi\pm L\rho$ as
in Lemma~\ref{bdy:barriers} bounds $v_n(0)$. The same boundary
identity then gives the upper bound.
\end{proof}

\subsection{A weighted integral estimate}
\label{wp:weighted}

Choose a collar $0\le r<2r_0$ where $d=r-Jr^2$ and $|Dd|\ge1/2$,
with $r_0$ fixed by the boundary geometry. Take a smooth $0\le\zeta\le1$
supported there and equal to one on $0\le r\le r_0$. Extend $\nu$ to this collar by $\nu=Dr$ and set
$\Pi=I-\nu\otimes\nu$. For a symmetric matrix $A$, write
$A'=A|_{\nu^\perp}$ in an orthonormal tangential frame. For an
$(n-1)\times(n-1)$ symmetric matrix $B$, let $T(B)$ be its
Newton tensor of order $k-2$. For $A\in\Ga$ we use
\begin{equation}\label{eq:wp-tensors}
 T(B)=\bigl(\sigma_{k-1}^{\alpha\beta}(B)\bigr)_{\alpha,\beta<n},
 \qquad Q(A)=I-\bigl(I+(\tr A)T(A')\bigr)^{-1}.
\end{equation}
For any unit tangent vector $e$, choose the frame so that
$e=e_1$ and $\nu=e_n$. Differentiation in the first diagonal entry gives
\begin{equation}\label{eq:wp-minor}
 e\cdot T(A')e=\sigma_{k-2}\bigl((A_{ij})_{2\le i,j<n}\bigr).
\end{equation}
If $H\in\Gamma_j$ and $|\xi|=1$, then, for $1\le i<j$,
\[
 \sigma_i(H|_{\xi^\perp})=\sigma_{i+1}^{ab}(H)\xi_a\xi_b>0
\]
by ellipticity of $\sigma_{i+1}$ on $\Gamma_{i+1}$. Thus
$H|_{\xi^\perp}\in\Gamma_{j-1}$ for $j\ge2$.
In particular $T(A')>0$, $0<Q(A)\le I$, and
$Q(A)\le (\tr A)T(A')$. By \eqref{eq:wp-minor} and the concavity
and homogeneity of $\sigma_{k-2}^{1/(k-2)}$,
$(e\cdot T(A')e)^{1/(k-2)}$ is superadditive in $A\in\Ga$.
In particular $T(A')$ increases as a quadratic form upon addition
of a matrix in $\overline{\Ga}$ to $A$.

Let $u,v$ be smooth on $\overline\Om$, admissible in $\Om$, with common boundary value
$\varphi$. Set $A=D^2u$, $B=D^2v$, $g=G(A)$, $h=G(B)$, $z=v-u$
and $d\mu=d^{1/2}\,dx$. Define
\[
\cD=\sigma_k^{ij}(A)B_{ij}-kg^{k-1}h
     +\sigma_k^{ij}(B)A_{ij}-kh^{k-1}g .
\]
Concavity and homogeneity of $G$ give
$\sigma_k^{ij}(A)B_{ij}\ge kg^{k-1}h$ and $\sigma_k^{ij}(B)A_{ij}\ge kh^{k-1}g$.
For the quadratic comparison identity see Chen, Zhou and
Zhu~\cite[Proposition~3.8]{CZZ}, and for its boundary-weighted form
see~\cite[Section~2]{LyuK2}.

\begin{lemma}\label{lem:weighted}
Suppose $m\le g,h\le M$. Then
\begin{equation}\label{eq:weighted-endpoints}
\int_\Om\cD\dm+
c\int_\Om\frac{\zeta}{r^2}\Pi Dz\cdot\bigl(T(A')+T(B')\bigr)\Pi Dz\dm
\le C\int_\Om|g-h|^2\dm .
\end{equation}
On $\{\cD\le1\}$, there is a constant $C=C(n,k,m,M)$ such that
\begin{equation}\label{eq:cone-comparison}
 \begin{aligned}
& C^{-1}\tr A\le \tr B\le C\tr A,\\
& C^{-1}T(A')\le T(B')\le CT(A'),\\
& C^{-1}Q(A)\le Q(B)\le CQ(A).
 \end{aligned}
\end{equation}
The matrix inequalities hold in the collar.
\end{lemma}

\begin{proof}
For $s>0$ put $\psi_s=(d+s)^{1/2}-s^{1/2}$,
$A_t=(1-t)A+tB$ and $J_s(t)=\int_\Om\psi_s\sigma_k(A_t)\,dx$.
The identity $\partial_i\sigma_k^{ij}(D^2w)=0$ follows by expanding
the principal minors and commuting third derivatives
(see Reilly~\cite[Proposition~2.1]{Reilly}). Since $z=\psi_s=0$ on
$\partial\Om$, two integrations by parts give
\[
 J_s'(t)=\int_\Om\psi_s\sigma_k^{ij}(A_t)z_{ij}\,dx
        =\int_\Om z\sigma_k^{ij}(A_t)(\psi_s)_{ij}\,dx .
\]
Differentiating $\partial_i\sigma_k^{ij}(A_t+\tau D^2\psi_s)=0$ at
$\tau=0$ gives
$\partial_i\bigl(\sigma_k^{ij,\ell m}(A_t)(\psi_s)_{\ell m}\bigr)=0$.
Thus differentiation in $t$ and integration by parts, using $z=0$
on the boundary and the symmetry in the two pairs of indices, yield
\begin{equation}\label{eq:polarized-secondvariation}
 J_s''(t)=-\int_\Om\sigma_k^{ij,\ell m}(A_t)
                      (\psi_s)_{\ell m}z_i z_j\,dx .
\end{equation}

Now
\[
 -D^2\psi_s=
 \tfrac12(d+s)^{-1/2}D^2\rho+
 \tfrac14(d+s)^{-3/2}Dd\otimes Dd .
\]
For a unit vector $\xi$ and $H,E\in\Ga$, choose an orthonormal
frame with last vector $\xi$. Differentiation in the normal entry gives
\[
 \sigma_k^{ij,rs}(H)E_{ij}\xi_r\xi_s
 =\sigma_{k-1}^{\alpha\beta}(H|_{\xi^\perp})E_{\alpha\beta}\ge0 .
\]
The inequality follows from concavity and homogeneity of
$\sigma_{k-1}^{1/(k-1)}$ on $\Gamma_{k-1}$ and extends by continuity
to $E\in\overline{\Ga}$. Thus both terms in $-D^2\psi_s$ have
nonnegative contributions throughout $\Om$. In the collar the same
identity in a frame with last vector $\nu$ gives
\[
 \sigma_k^{ij,rs}(H)\nu_i\nu_j\xi_r\xi_s=\Pi\xi\cdot T(H')\Pi\xi .
\]
By \eqref{eq:wp-minor} and superadditivity,
\[
 T\bigl((1-t)A'+tB'\bigr)\ge(1-t)^{k-2}T(A')+t^{k-2}T(B') .
\]
On the other hand, direct expansion of $J_s'(1)-J_s'(0)$ yields
\[
 \int_\Om\psi_s\cD\,dx+
 \int_0^1J_s''(t)\,dt
 =k\int_\Om\psi_s(g-h)(g^{k-1}-h^{k-1})\,dx .
\]
Retain the term with $Dd\otimes Dd$ on the support of $\zeta$, where
$Dd$ is parallel to $\nu$ and $|Dd|\ge1/2$. Integrating in $t$ and
letting $s\downarrow0$, we obtain \eqref{eq:weighted-endpoints} by
monotone convergence on the left and dominated convergence on the right,
since $d\asymp r$.

On $\{\cD\le1\}$ there is a $C_1=C_1(n,k,m,M)$ such that
$C_1A-B$ and $C_1B-A$ belong to $\overline{\Ga}$. Indeed, by G{\aa}rding's theorem~\cite{Garding}, the numbers $\lambda_j$
in $\sigma_k(A+tB)=\sigma_k(A)\prod_{j=1}^k(1+t\lambda_j)$ are positive,
with sum $\sigma_k^{ij}(A)B_{ij}/\sigma_k(A)$. Since each of the two differences
in $\cD$ is nonnegative, $\cD\le1$ bounds this sum by $C_1$. The
largest $\lambda_j$ is therefore at most $C_1$, which gives
$C_1A-B\in\overline{\Ga}$. Interchanging $A$ and $B$ gives the other
inclusion. Taking traces gives the first line of
\eqref{eq:cone-comparison}. By \eqref{eq:wp-minor},
$T(B')\le C_1^{k-2}T(A')$, and the reverse inequality follows by
interchanging $A$ and $B$. Finally,
$H\mapsto I-(I+H)^{-1}$ is increasing on positive definite matrices.
For $a\ge1$, the inequality $ax/(1+ax)\le ax/(1+x)$ for $x>0$
gives the stated comparison of $Q(A)$ and $Q(B)$.
\end{proof}

\subsection{An auxiliary Dirichlet problem}
\label{wp:auxiliary}

For $p\in\R^n$ and $A\in\Ga$, put $y=\Pi(p-D\Phi)$ and define
\[
Z(x,p,A)=\frac{\zeta\, y\cdot Q(A)y}{r^2},\qquad
Z_\eps(x,p,A)=\frac{\zeta\, y\cdot Q(A)y}{(r+\eps)^2}.
\]
Set $Z=Z_\eps=0$ off the collar and write $Z_v=Z(x,Dv,D^2v)$ for
admissible $v$. If $v$ is smooth with boundary value $\varphi$, then
$\Pi D(v-\Phi)=0$ on $\partial\Om$ and $Z_v$ is bounded.
Fix $P\ge1$ and $\eps>0$ smaller than half the collar width. Choose a
smooth nondecreasing function $\eta$ equal to zero on $(-\infty,1]$ and
one on $[2,\infty)$, and a smooth cutoff $0\le\chi_\eps\le1$ equal
to one where $r\le\eps$ and zero where $r\ge2\eps$.
For smooth functions $m\le b\le g\le M$, set
\begin{equation}\label{eq:aux-density}
h_\eps(x,p,A)=g-(g-b)\left[
\chi_\eps+(1-\chi_\eps)\eta\left(\frac{\tr A+Z_\eps(x,p,A)}{P}\right)
\right].
\end{equation}
Consider the Dirichlet problem
\begin{equation}\label{eq:aux-problem}
G(D^2v)=h_\eps(x,Dv,D^2v)\quad\text{in }\Om,\qquad
v=\varphi\quad\text{on }\partial\Om .
\end{equation}
For every admissible solution, $G(D^2v)=b$ on
$\{\Delta v>2P\}\cup\{Z_v>8P\}$. Indeed, $\chi_\eps=1$ if $r\le\eps$,
and $Z_\eps(x,Dv,D^2v)\ge Z_v/4$ otherwise.

\begin{proposition}\label{prop:wp-aux}
Suppose the domain and boundary values are smooth, and use a smooth
extension $\Phi$ satisfying the bounds above. There is a constant
$\delta_*>0$, depending only on $n,k,m,M$ and $\eta$, such that
\eqref{eq:aux-problem} has a smooth admissible solution whenever
$\|g-b\|_\infty\le\delta_*$. Every such solution satisfies
\begin{equation}\label{eq:auxbound}
\sup_\Om(\Delta v+Z_v)\le CP+C_b.
\end{equation}
In a fixed collar of $\partial\Om$ one has
$|\Pi D(v-\Phi)|/r\le C\sqrt P+C_b$.
\end{proposition}

\begin{proof}
\emph{The a priori estimate.}
At a point of the collar, let $A\in\Ga$ satisfy $m\le G(A)\le M$ and set $B=A'$.
If $\tau$ is a unit vector
in $\nu^\perp$ and $|B\tau|\le C_1$, then
\begin{equation}\label{eq:cap-direction}
 \tau\cdot Q(A)\tau\ge c(n,k,m,C_1)>0 .
\end{equation}
The expansion by the normal row and column gives
\[
 m^k\le \sigma_k(A)\le\sigma_k(B)+A_{\nu\nu}\sigma_{k-1}(B)
 \le C(\tr A)\sigma_{k-1}(B).
\]
Indeed, if $\sigma_k(B)>0$, then $B\in\Gamma_k$ and the Maclaurin
inequalities give $\sigma_k(B)\le C(\tr B)\sigma_{k-1}(B)$.
Also $|A_{\nu\nu}|+|\tr B|\le C\tr A$ since $|A|\le\tr A$.
In an eigenbasis of $B$, let $\lambda_j$ and
$\mu_j=\sigma_{k-2}(\lambda|j)>0$ be the corresponding eigenvalues of
$B$ and $T(B)$, respectively, where
$\lambda|j$ denotes omission of $\lambda_j$. Expanding at that entry gives
\[
 \sigma_{k-1}(B)\le |\lambda_j|\mu_j+C\mu_j^{(k-1)/(k-2)}.
\]
Here it suffices to use the Maclaurin inequality when
$\sigma_{k-1}(\lambda|j)>0$, in which case $\lambda|j\in\Gamma_{k-1}$.
If $\mu_j\le1$, the two displayed inequalities give
$\bigl((\tr A)\mu_j\bigr)^{-1}\le C(1+|\lambda_j|)$, and this also holds for
$\mu_j>1$ since $\tr A\ge c(n,k)m$. Thus
\[
 \tau\cdot Q(A)^{-1}\tau
 =\sum_j\bigl(1+\bigl((\tr A)\mu_j\bigr)^{-1}\bigr)|\tau_j|^2
 \le C(1+|B\tau|).
\]
Since $1\le(\tau\cdot Q(A)\tau)(\tau\cdot Q(A)^{-1}\tau)$,
this proves \eqref{eq:cap-direction}.

On the collar $0\le r\le r_0$, set $Y=\Pi D(v-\Phi)$ and $\psi(y)=|y|+|y|^2$.
Lemma~\ref{bdy:barriers} gives $|D(v-\Phi)|\le C$. We claim that
\begin{equation}\label{eq:cap-vector-barrier}
 \psi(Y)\le H d,\qquad H=C(1+\sqrt P+\|Db\|_\infty).
\end{equation}
Since $Y=0$ on $\partial\Om$ and $d\asymp r$, the inequality holds
on both edges of the collar when $C$ is sufficiently large.
At a positive interior maximum of $\psi(Y)-Hd$ we have $Y\ne0$.
Write $\tau=Y/|Y|$ and $A=D^2v$ there.
For each unit tangent vector $e$ to the level surface of $r$, the
critical point condition and $D_ed=0$ give
\[
 0=\tau\cdot D_eY
 =e\cdot D^2(v-\Phi)\tau
  +\tau\cdot(D_e\Pi)D(v-\Phi) .
\]
It follows that $|A'\tau|\le C$.
Since $\psi(Y)\le C|Y|$ and $d\ge c r$, \eqref{eq:cap-direction}
gives $Z_v\ge c|Y|^2/r^2>8P$ for the stated choice of $H$ with $C$
sufficiently large. Hence $G(D^2v)=b$ in a neighborhood of the maximum.

Write $a^{ij}=G^{ij}(A)$, $L=a^{ij}\partial_{ij}$, and diagonalize $A$ at
this point in orthonormal coordinates $e_i$, with eigenvalues
$\lambda_i$ and $a_i$ for $A$ and $a$, respectively.
Since the distance normals are constant along their normal rays,
$D_\nu\Pi=0$ and $|D_{e_i}\Pi|\le C|\Pi e_i|$.
Differentiating $G(D^2v)=b$ gives $L(Dv)=Db$. Applying the
product rule to $Y=\Pi D(v-\Phi)$, we obtain
\[
 |LY|\le \|Db\|_\infty+C\tr a
              +C\sum_i a_i|\lambda_i|\,|\Pi e_i|,
 \qquad
 |D_{e_i}Y|^2\ge \tfrac12\lambda_i^2|\Pi e_i|^2-C .
\]
The second inequality holds for each $i$. The first uses the
bounded first two derivatives of $\Pi$
and the bounded first three derivatives of $\Phi$.
At $Y\ne0$ one has $|D\psi(Y)|\le C$ and $D^2\psi(Y)\ge 2I$.
The chain rule and Young's inequality thus imply
\[
 L(\psi(Y))\ge -C\tr a-C\|Db\|_\infty.
\]
By \eqref{eq:linear-barrier}, $Ld\le -\theta\tr a$ and
$\tr a\ge G(I)>0$. Thus $L(\psi(Y)-Hd)>0$ at the maximum,
a contradiction. This proves \eqref{eq:cap-vector-barrier}.
It gives $Z_v\le CP+C_b$ near the boundary. On the remaining support of $\zeta$, use $Q(A)\le I$,
the gradient bound and $r\ge r_0$. On the boundary, $D_\nu\Pi=0$,
and division of $|Y|\le Hd$ by $r$ along the normal yields
$|(D^2v)_{T\nu}|\le C\sqrt P+C_b$.

The tangential block $B=(D^2v)|_{\nu^\perp}$ is bounded on the
boundary and satisfies $\sigma_{k-1}(B)\ge c_0$ by Proposition~\ref{bdy:margin}.
For $z=(D^2v)_{T\nu}$ the block expansion is
\[
 \sigma_{k-1}(B)v_{\nu\nu}
  =G(D^2v)^k-\sigma_k(B)+\sigma_{k-1}^{\alpha\beta}(B)z_\alpha z_\beta.
\]
Thus $\Delta v\le CP+C_b$ on $\partial\Om$. Choose $C_b$ so that
$2C_bG(I)>\|\Delta b\|_\infty$. Where $\Delta v>2P$, we have
$G(D^2v)=b$, and concavity gives $L(\Delta v)\ge\Delta b$.
Fix $x_*\in\Om$. Then $L(\Delta v+C_b|x-x_*|^2)>0$ there. The maximum principle and
the boundary estimate give \eqref{eq:auxbound}.

\emph{Solvability.}
Write $F(x,p,A)=\log\sigma_k(A)-k\log h_\eps(x,p,A)$, so that
\eqref{eq:aux-problem} is $F(x,Dv,D^2v)=0$. Fix $x,p$ and let $H$
be symmetric. By hyperbolicity we may write
\[
 \begin{aligned}
 \sigma_k(A+tH)&=\sigma_k(A)\prod_{j=1}^k(1+t\lambda_j),\\
 N^2=\sum_j\lambda_j^2
     &=-\left.\frac{d^2}{dt^2}\log\sigma_k(A+tH)\right|_{t=0}.
 \end{aligned}
\]
If $N=0$, then $\sigma_k(A+tH)=\sigma_k(A)>0$ for all $t$, so
$A+tH$ stays in $\Ga$. Positivity of $\tr(A+tH)$ for all $t$ gives
$\tr H=0$, and $|A+tH|\le\tr A$ for all $t$ gives $H=0$. Otherwise normalize $H$ so that $N=1$.
Then $A+tH\in\Ga$ for $|t|<1$, and $A\pm H\in\overline\Ga$.

In the collar put $V(t)=\tr(A+tH)\,T(A'+tH')$. For $|t|\le1$ both
$A+tH$ and $2A-(A+tH)$ belong to $\overline\Ga$. The trace and the
matrix $T(A')$ both increase upon addition of a matrix from
$\overline\Ga$, by \eqref{eq:wp-minor}. Homogeneity gives
$0\le V(t)\le 2^{k-1}V(0)$. For each tangent vector $e$, the polynomial
$e\cdot V(t)e$ has degree at most $k-1$.
Its first two derivatives at zero are bounded by $C(k)$ times
its supremum on $[-1,1]$, by equivalence of norms on polynomials of
bounded degree. Thus $-CV\le V^{(j)}\le CV$, $j=1,2$, at $t=0$.

Set $R(t)=(I+V(t))^{-1}$. Since $Q(A+tH)=I-R(t)$,
differentiation at $t=0$ gives
\[
 \left.\frac{d}{dt}Q(A+tH)\right|_{t=0}=RV'R,\qquad
 \left.\frac{d^2}{dt^2}Q(A+tH)\right|_{t=0}
       =RV''R-2RV'RV'R ,
\]
where the matrices on the right are evaluated at zero.
Since $R\le V^{-1}$ and $-CV\le V'\le CV$, we have
$0\le V'RV'\le V'V^{-1}V'\le C^2V$. Also $RVR\le Q(A)$.
Thus the first two derivatives of $Q(A+tH)$ at zero lie between
$-CQ(A)$ and $CQ(A)$ as quadratic forms. The normalization also
gives $|\tr H|\le\tr A$. Restoring $N$ and writing
$z=\tr A+Z_\eps(x,p,A)$, we obtain $|z'|\le C Nz$ and
$|z''|\le C N^2z$ at $t=0$ along $A+tH$. Off the collar the same
bounds hold since $z=\tr A$.
Only $P\le z\le2P$ contributes when $\eta(z/P)$ is differentiated.
Since $h_\eps\ge m$, formula \eqref{eq:aux-density} therefore gives
\[
 |(\log h_\eps)'|\le C\delta_*N,\qquad
 |(\log h_\eps)''|\le C\delta_*N^2
\]
when $\|g-b\|_\infty\le\delta_*\le1$, independently of $P$ and $\eps$.
For sufficiently small $\delta_*$ this gives
$F^{ij,rs}H_{ij}H_{rs}\le-N^2/2$ for every symmetric $H$. If $H\ge0$ the
relative eigenvalues are nonnegative, so that
$\frac{\sigma_k^{ij}(A)H_{ij}}{\sigma_k(A)}=\sum_j\lambda_j\ge N$ and
\[
 F^{ij}H_{ij}\ge\frac{\sigma_k^{ij}(A)H_{ij}}{2\sigma_k(A)}.
\]

For the method of continuity replace $g$ by $b+s(g-b)$, $0\le s\le1$.
At $s=0$ the smooth admissible solution with $G(D^2v)=b$ exists
by \cite[Theorem~3]{CNS}. The linearized Dirichlet operator is elliptic with no
zeroth-order term, so the Schauder theory and the implicit function
theorem give openness. For fixed $P,\eps$ and smooth data,
\eqref{eq:auxbound} and $G(D^2v)\ge m$ confine the Hessians to a
compact subset of $\Ga$ along the path. There $F$ is uniformly elliptic and concave in $A$. By Lemma~\ref{bdy:barriers}, the gradients are bounded as well.
For fixed $P,\eps$, the first and second derivatives of $F$ are thus
bounded on a common compact neighborhood of the values of $(x,Dv,D^2v)$
along the path. Krylov's global estimate~\cite{KrylovBoundary}, in the form
of Gilbarg and Trudinger~\cite[Theorem~17.26$'$]{GT}, gives a uniform $C^{2,\beta}(\overline\Om)$
bound for some $\beta>0$. Compactness in $C^{2,\beta/2}(\overline\Om)$
gives closedness, and Schauder estimates give smoothness. These last
bounds may depend on $P,\eps$ and the fixed smooth data.
\end{proof}

\subsection{Estimates for level sets}
\label{wp:tails}

\begin{proof}[Proof of Theorem~\ref{thm:wp}]
First assume that the domain, data and solution are smooth.
We use the level set argument of~\cite{LyuInterior,LyuK2}.
Fix $0<\delta\le\delta_*$ and take a smooth function $b$ such that
$m\le b\le g\le M$ and $g-b\le\delta$. For $P\ge1$ and $\eps>0$
sufficiently small let $v$ solve \eqref{eq:aux-problem}, and put
$A=D^2u$, $B=D^2v$. For $w=u,v$ write
\[
 y_w=\frac{\sqrt\zeta}{r}\Pi D(w-\Phi),\qquad
 E_w(t)=\{\Delta w+Z_w>t\},
\]
with $y_w=0$ off the collar. In the collar $Z_u=y_u\cdot Q(A)y_u$ and
$Z_v=y_v\cdot Q(B)y_v$. By \eqref{eq:auxbound}, $\tr B+Z_v\le LP$ when
$P\ge P_b$, where $L$ depends only on the data fixed at the beginning
of the section and $P_b$ may also depend on $\|b\|_{C^2}$.
Since $Z_\eps(x,Dv,B)\le Z_v$, the right-hand side of \eqref{eq:aux-problem}
agrees with $g$ outside $E_v(P)\cup\{r<2\eps\}$. Lemma~\ref{lem:weighted}
gives
\begin{equation}\label{eq:energy-source}
 \int_\Om\!\bigl(\cD+(y_u-y_v)\cdot\bigl(T(A')+T(B')\bigr)(y_u-y_v)\bigr)\dm
 \le C\delta^2\bigl(\mu(E_v(P))+\mu\{r<2\eps\}\bigr),
\end{equation}
where the quadratic form in the integrand is taken to be zero off the collar.

Let $F$ be the set where the integrand on the left of
\eqref{eq:energy-source} exceeds one. On $\Om\setminus F$ we have
$C^{-1}\tr A\le\tr B\le C\tr A$ by \eqref{eq:cone-comparison}.
In the collar, $Q(A)\le CQ(B)$ and $Q(A)\le(\tr A)T(A')$, so
\[
 \begin{aligned}
 Z_u&\le 2y_v\cdot Q(A)y_v
       +2(y_u-y_v)\cdot Q(A)(y_u-y_v)\\
    &\le C Z_v+2\tr A .
 \end{aligned}
\]
Off the strip both $Z_u$ and $Z_v$ vanish. Interchanging $u$ and $v$
therefore gives, with a uniform $C_0\ge1$,
\begin{equation}\label{eq:E-control}
 C_0^{-1}(\Delta v+Z_v)\le\Delta u+Z_u\le C_0(\Delta v+Z_v)
 \qquad\text{on }\Om\setminus F .
\end{equation}
Hence $E_v(P)\subset E_u(P/C_0)\cup F$, while the uniform bound
$\Delta v+Z_v\le LP$ gives $E_u(C_0LP)\subset F$. Substitution in
\eqref{eq:energy-source} gives
\[
 \mu(F)\le C\delta^2\bigl(\mu(E_u(P/C_0))+\mu(F)+\mu\{r<2\eps\}\bigr).
\]
Choose $\delta$ small enough that $C\delta^2\le1/2$. Absorb
$C\delta^2\mu(F)$ and let $\eps\downarrow0$ at fixed
$P\ge P_b$. Since $\mu\{r<2\eps\}\to0$ and $E_u$ is independent of
$\eps$, we obtain
\[
 \mu(E_u(C_0LP))\le C\delta^2\mu(E_u(P/C_0)).
\]
Here $C,C_0,L$ are independent of $P,\eps$ and the derivatives of $g,b$.

Let $q=2p$. Integration in $P\ge P_b$ against $qP^{q-1}\,dP$ gives
\[
 \int_\Om (\Delta u+Z_u)^q\dm
 \le (C_0LP_b)^q\mu(\Om)+C\delta^2(C_0^2L)^q
             \int_\Om (\Delta u+Z_u)^q\dm .
\]
Both integrals are finite for smooth $u$, since $Z_u$ is bounded.
Choose $\delta$ so that the last coefficient is at most $1/2$.
With $b$ and $P_b$ chosen for this $\delta$, the resulting bound and
Cauchy--Schwarz give
\[
 \int_\Om (\Delta u)^p\,dx
 \le \left(\int_\Om (\Delta u)^{2p}\dm\right)^{1/2}
      \left(\int_\Om d^{-1/2}\,dx\right)^{1/2}.
\]
The second factor is finite since $d\asymp r$. The inclusion
$\Ga\subset\Gamma_2$ gives $|D^2u|\le\Delta u$, and
Lemma~\ref{bdy:barriers} bounds $u$ and $Du$.
This proves the estimate for smooth solutions.

We pass to continuous data. Fix $1<p<\infty$, set $m=\lambda^{1/k}/2$, $M=2\Lambda^{1/k}$,
and take $\delta$ as above, reducing it so that $\delta\le m/2$.
Extend $g=f^{1/k}$ to a fixed neighborhood of $\overline\Om$ using
a concave majorant of its modulus of continuity, and truncate the
extension to $[\lambda^{1/k},\Lambda^{1/k}]$. Mollifying to within
$\delta/8$ of $g$ and subtracting $\delta/2$ gives a single smooth
function $b\ge m$ with $3\delta/8\le g-b\le5\delta/8$. Its $C^2$ norm
is bounded in terms of $\delta,m,M,\Om$ and the modulus of continuity
of~$g$. Thus $P_b$ is fixed throughout the approximation.

Approximate the defining function and the extension $\Phi$ in $C^3$
by smooth functions. By subtracting a constant tending to zero from
each approximating defining function, we may take the smooth domains
$\Om_j$ to contain $\overline\Om$. Their boundaries converge in $C^3$
in fixed graph coordinates. Their principal curvature vectors lie in a
fixed compact subset of $\Gamma_{k-1}$. The third derivatives of the
graphs and of the approximations $\Phi_j$ have a common modulus of
continuity.
Let $g_j$ be smooth mollifications of the extended $g$ converging
uniformly to $g$. Then $b\le g_j\le b+\delta$ for large $j$.
By Caffarelli, Nirenberg and Spruck~\cite[Theorem~3]{CNS}, there are smooth
admissible solutions $u_j$ with density $g_j^k$ in $\Om_j$ and boundary
values $\Phi_j$. The preceding estimate bounds their restrictions to
$\Om$ in $W^{2,p}(\Om)$.

By weak compactness and the gradient estimate in Lemma~\ref{bdy:barriers},
a subsequence converges weakly in $W^{2,p}(\Om)$ and locally uniformly
in $\Om$ to a $k$-convex function $u_*$. The boundary barriers in
the same lemma give
$|u_j-\Phi_j|\le C\dist(\cdot,\partial\Om_j)$, with $C$ independent
of $j$. Passing to the limit shows that $u_*$ extends continuously to
$\partial\Om$ with value $\varphi$. Weak continuity of the Hessian
measure under locally uniform convergence~\cite[Theorem~1.1(ii)]{TW1}
gives $\mu_k[u_*]=f\,dx$. By
comparison~\cite[Theorem~3.1]{TW1}, $u_*=u$. Weak lower semicontinuity
gives \eqref{eq:main-estimate}.
\end{proof}

\section{Global \texorpdfstring{$C^{2,\alpha}$}{C2,alpha} estimates}\label{sec:c2alpha}

Write $g=f^{1/k}$, $0<m\le g\le m_1$, and $S=\Delta u$.
Constants may depend on the data in Theorem~\ref{thm:c2alpha},
and $C_\eps$ may also depend on $\eps$.

\subsection{An interior estimate for the Laplacian}\label{c2:persistence}

Let $R_s=(I-s\Delta)^{-1}$ on $\R^n$, as in
Chen, Zhou and Zhu~\cite[Section~2.1]{CZZ} and~\cite{LyuK2}.
The Hessian quotient
\begin{equation}\label{c2:eq:quotient}
 A\longmapsto\left(\frac{\sigma_k(A)}{\tr A}\right)^{1/(k-1)}
 \qquad(A\in\Ga)
\end{equation}
is elliptic and concave on $\Ga$~\cite{Trudinger1990},
as stated in~\cite[(1.10), (1.11) and (1.15)]{Trudinger1995}.

\begin{lemma}\label{c2:lem:persistence}
Let $u$ be smooth and admissible on $\overline{B_R(x_0)}$, $0<R\le1$,
with $G(D^2u)=g$. Suppose
\[
 S\le C_1M\quad\text{in }B_R(x_0),\qquad
 S(x_0)\ge M-e,\qquad M\ge1,\quad 0\le e\le M/8.
\]
There are positive constants $c,C,\eta_0$ depending only on
$n,k,\alpha,m,m_1,[g]_{C^\alpha}$ and $C_1$ with the following
property. If $0<h\le cR$ and $h^{\alpha/2}+M^{k/2}h/R\le\eta_0$,
there is a nonnegative integrable function $K$, compactly supported
in $B_R(x_0)$, with $\int K\le1$, such that
\begin{equation}\label{c2:eq:trace-average}
 \int K\Delta u\ge M-e-CM\bigl(h^{\alpha/2}+M^{k/2}h/R\bigr).
\end{equation}
Moreover, for $p>n/2$ and $w\in C^2(\overline{B_R(x_0)})$,
\[
 \|K\|_{L^{p/(p-1)}(\R^n)}\le C_p h^{-n/p},\qquad
 \left|\int K\Delta w\right|\le Ch^{-2}\|w\|_{L^\infty(B_R(x_0))}.
\]
\end{lemma}

\begin{proof}
Choose a smooth cutoff $0\le\chi\le1$ supported in $B_R(x_0)$,
equal to one on $B_{3R/4}(x_0)$, with $|D^j\chi|\le CR^{-j}$ for
$j=1,2$. Set $\widetilde u=\chi(u-\ell)$, where $\ell$ is the affine
tangent function of $u$ at $x_0$. Since $|D^2u|\le S$, Taylor's
formula gives $|D^2\widetilde u|\le CM$. The kernel of $R_s$ is
\[
 B_s(y)=\int_0^\infty e^{-t}(4\pi st)^{-n/2}e^{-|y|^2/(4st)}\,dt.
\]
It has mass one, and for $r=\sqrt{s}$ its moment of order $q>0$ is
$C_{n,q}r^q$. In particular, its mass outside $B_{R/2}(0)$ is at most $C(r/R)^2$.

Set $S_s=R_s\Delta\widetilde u$ and
\[
 N(s)=\max_{x\in\R^n}
      \left(S_s(x)-AM|x-x_0|^2/R^2\right).
\]
Since $|S_s|\le CM$, choosing $A$ sufficiently large places every
maximizer $x_s$ in $B_{R/4}(x_0)$. At $x_s$ we have
$D^2S_s\le2AMR^{-2}I$. Applying $(I-s\Delta)R_s=I$ to
$D^2\widetilde u$ gives
\[
 R_sD^2\widetilde u(x_s)
 =D^2u(x_s)+sD^2S_s(x_s)\le D^2u(x_s)+CM(r/R)^2I .
\]

Suppose $N(s)\ge M/2$ and $r\le cR$. Normalize
$B_s(x_s-y)\,dy$ on $B_{3R/4}(x_0)$ to obtain a probability measure
$\mu$. Set $\bar H=\int D^2u\,d\mu$ and $\bar S=\tr\bar H$.
The tail estimate gives
$|\bar H-R_sD^2\widetilde u(x_s)|\le CM(r/R)^2$, so
$\bar S\asymp M$ and $\bar H\le D^2u(x_s)+CM(r/R)^2I$.
The Hessians on $\overline{B_{3R/4}(x_0)}$ form a compact subset of
the convex cone $\Ga$, so $\bar H\in\Ga$. By monotonicity of $\sigma_k$
and the bounds $|\sigma_k^{ij}|\le CM^{k-1}$ for all $i,j$ on the segment from $D^2u(x_s)$
to $D^2u(x_s)+CM(r/R)^2I$, we obtain
\[
 \sigma_k(\bar H)\le g(x_s)^k+CM^k(r/R)^2.
\]

Put $b=1/(k-1)$, $\tau=S/\bar S$ and $v=g/g(x_s)$.
Then $0<\tau\le C$, $c\le v\le C$, and $\int \tau\,d\mu=1$.
The concavity in \eqref{c2:eq:quotient} gives
\begin{equation}\label{c2:eq:jensen}
 \int v^{1+b}\tau^{-b}\,d\mu
 \le \left(\frac{\sigma_k(\bar H)}{g(x_s)^k}\right)^b
 \le 1+CM^k(r/R)^2.
\end{equation}
The second derivative of $x^{-b}$ has a positive lower bound on
the range $0<x\le C$ of $\tau/v$. Taylor's formula at $1$, multiplied
by $v$, thus gives
\[
 v^{1+b}\tau^{-b}\ge (1+b)v-b\tau+c(\tau-v)^2 .
\]
By the H\"older bound and the kernel moments,
$\int |v-1|\,d\mu\le Cr^\alpha$ and
$\int |v-1|^2\,d\mu\le Cr^{2\alpha}$. Integration in the preceding
inequality and \eqref{c2:eq:jensen} therefore give
$\int(\tau-1)^2\,d\mu\le C(r^\alpha+M^k(r/R)^2)$.
On the complement of $B_{3R/4}(x_0)$ we use
$|\Delta\widetilde u-\bar S|\le CM$ and the kernel tail bound. Since the mean minimizes
the mean square
deviation, it follows that
\begin{equation}\label{c2:eq:variance}
 R_s\bigl[(\Delta\widetilde u-S_s(x_s))^2\bigr](x_s)
 \le CM^2\bigl(r^\alpha+M^k(r/R)^2\bigr).
\end{equation}

For bounded smooth $w$ and any constant $a$, integration by parts in
the heat parameter gives
\[
 s\partial_sR_sw=\int_0^\infty (t-1)e^{-t}e^{st\Delta}(w-a)\,dt,
 \qquad |s\partial_sR_sw|^2\le R_s[(w-a)^2].
\]
The inequality follows by Cauchy--Schwarz since
$\int_0^\infty(t-1)^2e^{-t}\,dt=1$. The maximum $N$ is locally
Lipschitz for $s>0$, and at each point of differentiability
$N'(s)=\partial_s S_s(x_s)$. By \eqref{c2:eq:variance},
\[
 |sN'(s)|\le CM\bigl(r^{\alpha/2}+M^{k/2}r/R\bigr)
 \qquad\text{whenever }N(s)\ge M/2 .
\]
As $s\downarrow0$, $S_s\to\Delta\widetilde u$ uniformly and
$N(0)\ge M-e\ge7M/8$. Integrating with $ds/s=2\,dr/r$ shows that
$N(s)\ge N(0)-CM(h^{\alpha/2}+M^{k/2}h/R)$ up to $s=h^2$ or the
first time $N=M/2$. Taking $\eta_0$ sufficiently small excludes this
second possibility.

At a maximizer $z$ for $N(h^2)$ put $K(y)=\chi(y)B_{h^2}(z-y)$.
The difference $\Delta\widetilde u-\chi\Delta u$ vanishes on
$B_{3R/4}(x_0)$ and is bounded by $CM$. Hence
$|S_{h^2}(z)-\int K\Delta u|\le CM(h/R)^2$.
Since $S_{h^2}(z)\ge N(h^2)$, this proves \eqref{c2:eq:trace-average}.

The heat kernel formula and Minkowski's inequality give, for $p>n/2$,
\[
 \|B_{h^2}\|_{L^{p/(p-1)}}\le
 C_p h^{-n/p}\int_0^\infty e^{-t}t^{-n/(2p)}\,dt,\qquad
 \|DB_{h^2}\|_{L^1}\le Ch^{-1}.
\]
Also $h^2\Delta B_{h^2}=B_{h^2}-\delta_0$ in distributions. By
the product rule, $\Delta K$ is therefore a signed measure of total
variation at most $C(h^{-2}+R^{-1}h^{-1}+R^{-2})\le Ch^{-2}$.
Integration by parts gives the stated bound for $\int K\Delta w$.
\end{proof}

\begin{lemma}\label{c2:lem:subpower}
Under the hypotheses of Theorem~\ref{thm:c2alpha}, assume that the domain
and the data are smooth and that $u$ is smooth on $\overline\Om$. Then for every $\eps>0$,
\begin{equation}\label{c2:eq:subpower}
 |D^2u(x)|\le\Delta u(x)\le C_\eps\dist(x,\partial\Om)^{-\eps}
 \qquad(x\in\Om).
\end{equation}
\end{lemma}

\begin{proof}
Write $d_0(x)=\dist(x,\partial\Om)$ and maximize $d_0^\eps S$ at
$x_0$. With $M=S(x_0)$ and $R=\min\{1,d_0(x_0)/2\}$ we have
$S\le2^\eps M$ on $B_R(x_0)$. We may assume $M\ge1$.
Take $h=c_*RM^{-k/2}$ with $c_*>0$ sufficiently small. Applying
Lemma~\ref{c2:lem:persistence} with $e=0$, then H\"older's inequality
and Theorem~\ref{thm:wp}, gives
\[
 M/2\le\int K\Delta u\le C_ph^{-n/p}\|\Delta u\|_{L^p(\Om)}
       \le C_pR^{-n/p}M^{kn/(2p)}.
\]
For $p>kn/2+n/\eps$ this gives $d_0(x_0)^\eps M\le C_\eps$, as
claimed.
\end{proof}

\subsection{Smooth approximation}\label{c2:approximation}

Take the defining function $\rho$ of Lemma~\ref{bdy:barriers} and
the extension $\Phi$ fixed in Section~\ref{wp:nondegeneracy}.
Extend $\rho$ to a neighborhood
of $\overline\Om$ where $D^2\rho-\theta I\in\Ga$ for some $\theta>0$.
Let $\eta$ be a nonnegative even smooth kernel of mass one supported
in $B_{1/2}(0)$. With $\eta_t(x)=t^{-n}\eta(x/t)$, the kernel
$\zeta_t=(4\eta_t-\eta_{2t})/3$ has mass one and vanishing moments of orders
one through three. For $v=\rho,\Phi$ and $0\le j\le3$, Taylor's
formula gives
\begin{equation}\label{c2:eq:smoothing}
 \|D^j(v*\zeta_t-v)\|_\infty\le Ct^{3-j}\omega_3(t),\qquad
 \|D^4(v*\zeta_t)\|_\infty\le Ct^{-1}\omega_3(t),
\end{equation}
where $\omega_3$ is a common modulus of continuity for $D^3\rho,D^3\Phi$.

For $t=2^{-j}t_0$ set
\[
 \rho_t=\rho*\zeta_t-t^3,\qquad d_t=-\rho_t,\qquad \Phi_t=\Phi*\zeta_t.
\]
Let $\Om_t$ denote the domain bounded by $\rho_t=0$ near $\partial\Om$.
Choose $t_0$ small enough that $|\rho*\zeta_t-\rho|\le t^3/2$ for
$t\le2t_0$. Then $\rho_t<0$ on $\overline\Om$ and
$\rho_{2t}-\rho_t\le-5t^3/2$, so
\begin{equation}\label{c2:eq:nesting}
 \overline\Om\subset\Om_t,\qquad\overline{\Om_t}\subset\Om_{2t}.
\end{equation}
Moreover $D^2\rho_t-(\theta/2)I\in\Ga$. In a fixed boundary collar
$|Dd_t|\ge c$ and $d_t\asymp\dist(\cdot,\partial\Om_t)$. By the
implicit function theorem and \eqref{c2:eq:smoothing}, the local
boundary graphs change by $O(t^3)$ in the uniform norm and converge
in $C^2$. Thus the principal curvature vectors of $\partial\Om_t$
remain in a fixed compact subset of $\Gamma_{k-1}$.
The $C^3$ bounds and moduli of continuity of the third derivatives
of the graphs and of $\Phi_t$ are uniform, and their $C^4$ norms
are at most $Ct^{-1}$.

Extend $g$ to the same neighborhood preserving its bounds and H\"older
seminorm. We approximate the density separately near the boundary and
in the interior, as in~\cite[Section~2]{TWMA} and~\cite{LyuK2}.
Take a smooth function $0\le\chi\le1$ equal to one on
$(-\infty,2]$ and zero on $[3,\infty)$, and put
\[
 g_t=\chi(d_t/t)(g*\eta_t)+(1-\chi(d_t/t))(g*\eta_{t^{4/\alpha}}).
\]
Both convolutions have the same H\"older bound as $g$ and they differ
by at most $Ct^\alpha$. Thus the cutoff contributes at most
$Ct^\alpha\min\{1,|x-y|/t\}\le C|x-y|^\alpha$ to
$|g_t(x)-g_t(y)|$. In particular $m\le g_t\le m_1$ and
\begin{equation}\label{c2:eq:density}
 \begin{split}
 [g_t]_{C^\alpha(\overline{\Om_t})}&\le C,\qquad
 |D^jg_t|\le Ct^{\alpha-j}\quad(d_t<2t,\ j=1,2),\\
 |g_t-g|&\le Ct^\alpha\quad\text{in }\Om_t,\qquad
 |g_t-g|\le Ct^4\quad(d_t\ge3t).
 \end{split}
\end{equation}
The smooth Dirichlet theory~\cite[Theorem~3]{CNS} gives an admissible solution
$u_t$ of $G(D^2u_t)=g_t$ in $\Om_t$, with $u_t=\Phi_t$ on
$\partial\Om_t$.

Fix $\eps=\alpha/[100(k+2)]$ and $\beta=1-\eps/2$.
Taking $p>\max\{kn/2+n/\eps,2n/\eps\}$ in Theorem~\ref{thm:wp},
Sobolev embedding and Lemma~\ref{c2:lem:subpower} give
\begin{equation}\label{c2:eq:preparation}
 \|u_t\|_{C^{1,\beta}(\overline{\Om_t})}\le C_\eps,\qquad
 S_t:=\Delta u_t\le C_\eps\dist(\cdot,\partial\Om_t)^{-\eps}.
\end{equation}
Write $a^{ij}=G^{ij}(D^2u_t)$ and $\cL=a^{ij}\partial_{ij}$.
Logarithmic differentiation of \eqref{c2:eq:quotient} in any
nonnegative matrix direction gives
\begin{equation}\label{c2:eq:ellipticity}
 G^{ij}(A)\xi_i\xi_j\ge\frac{G(A)}{k\tr A}|\xi|^2 \qquad(A\in\Ga).
\end{equation}
Concavity and homogeneity of $G$ give $\tr a\ge G(I)$ and
$\cL d_t\le-(\theta/2)\tr a$. In the collar, \eqref{c2:eq:preparation}
and \eqref{c2:eq:ellipticity} give $a^{ij}(d_t)_i(d_t)_j\ge c_\eps d_t^\eps$.

We record the barrier used below. Suppose $d$ is positive in the
domain, vanishes on its boundary, and satisfies $\cL d\le0$ and
$a^{ij}d_i d_j\ge c_\eps d^\eps$ on $0<d<2s$. Fix a smooth cutoff
$0\le\chi_0\le1$ equal to one on $[0,1]$ and zero on $[2,\infty)$.
For $J\ge0$ let
\[
 \psi_{J,s}(d)=C_\eps J\int_0^{2s}\min\{d,y\}y^{-\eps}\chi_0(y/s)\,dy .
\]
This function is nondecreasing in $d$, vanishes at zero, and satisfies
$\psi_{J,s}''(d)=-C_\eps Jd^{-\eps}\chi_0(d/s)$ for $d>0$.
Taking $C_\eps$ sufficiently large gives
\begin{equation}\label{c2:eq:barrier}
 \begin{split}
 \cL\psi_{J,s}(d)&\le-J\chi_0(d/s),\qquad
 \psi_{J,s}'(0)\le C_\eps Js^{1-\eps},\\
 0\le\psi_{J,s}(d)&\le C_\eps Js^{1-\eps}\min\{d,s\}.
 \end{split}
\end{equation}

\begin{lemma}\label{c2:lem:preliminary}
For $M_t=\sup_{\Om_t}S_t$,
\begin{align}
 M_t&\le C_\eps t^{-\eps},\label{c2:eq:preliminary}\\
 M_t&\le\max\left\{\sup_{\partial\Om_t}S_t,\,
                       \sup_{\{d_t\ge t\}}S_t\right\}
                   +C_\eps t^{\alpha-\eps}.\label{c2:eq:localization}
\end{align}
\end{lemma}

\begin{proof}
By concavity, $\cL S_t\ge\Delta g_t\ge-Ct^{\alpha-2}$ on
$0<d_t<t$. Thus $\cL(S_t-\psi_{J,t}(d_t))\ge0$ there for
$J=Ct^{\alpha-2}$. Since $0\le\psi_{J,t}\le C_\eps t^{\alpha-\eps}$,
the maximum principle gives \eqref{c2:eq:localization}.
On $d_t\ge t$, the bound $|Dd_t|\le C$ and \eqref{c2:eq:preparation}
give $S_t\le C_\eps t^{-\eps}$. It remains to estimate $S_t$ on
$\partial\Om_t$.

At a boundary point take the inward normal as $e_n$ and any unit
tangent as $e_1$. Write the boundary as $x_n=b(x')$, with
$b(0)=D'b(0)=0$. Put $h=D'b_1(0)$ and
\begin{equation}\label{c2:eq:rigid}
 X=\partial_1+\sum_{i<n}h_i(x_i\partial_n-x_n\partial_i),
 \qquad w=u_t-\Phi_t .
\end{equation}
By the orthogonal invariance of $G$, differentiation along the rigid
motions generated by $X$ gives $\cL(Xu_t)=Xg_t$,
as in~\cite[(5.4)--(5.5)]{CNS}. Hence
$|\cL(Xw)|\le Ct^{\alpha-1}+C\tr a$ on $0<d_t<t$.
Since $w=0$ on the graph,
\[
 Xw=(h\cdot x'-b_1+b\,h\cdot Db)w_n=O(|x'|^2).
\]
By \eqref{c2:eq:preparation}, $|Xw|\le C_\eps(|x|^2+d_t^\beta)$
in a smaller boundary chart.

On $B_R\cap\{0<d_t<t\}$, with $R$ small and fixed, use
\[
 V=A d_t+B|x|^2+\psi_{J,t}(d_t),\qquad J=Ct^{\alpha-1}.
\]
On the graph, $|Xw|\le C|x'|^2$. On $\partial B_R$ the gradient
bound gives $|Xw|\le C$, and $R$ is fixed. Choose $B$ so that
$B|x|^2$ dominates both bounds. Then $A=C_\eps(1+t^{\beta-1})$, with a
sufficiently large constant, gives $V\ge|Xw|$ on $d_t=t$ and
$\cL V\le-|\cL(Xw)|$. Comparison yields $-V\le Xw\le V$, with
equality at the origin. There $\partial_n(Xw)=w_{1n}$ because $Dw$
is normal. Taking inward derivatives and using \eqref{c2:eq:barrier}
gives
\[
 |(u_t)_{1n}(0)|\le C_\eps(1+t^{\beta-1}+t^{\alpha-\eps})
                 \le C_\eps t^{-\eps/2}.
\]
For $B_t=((u_t)_{ij})_{i,j<n}$ and $z=((u_t)_{in})_{i<n}$,
expansion in the last row and column gives
\begin{equation}\label{c2:eq:block}
 (u_t)_{nn}\sigma_{k-1}(B_t)
 =g_t^k-\sigma_k(B_t)+\sigma_{k-1}^{\alpha\beta}(B_t)z_\alpha z_\beta.
\end{equation}
Since $|B_t|\le C$ and $\sigma_{k-1}(B_t)\ge c$ by
Proposition~\ref{bdy:margin}, we obtain
$S_t\le C(1+|z|^2)\le C_\eps t^{-\eps}$ on $\partial\Om_t$.
This proves \eqref{c2:eq:preliminary}.
\end{proof}

\subsection{Mixed second derivatives on the boundary}\label{c2:boundary}

Fix finitely many boundary charts in which $\partial\Om_t$ is given by
$x_n=b_t(x')$ for $x'\in B'_{3r}$, with $\Om_t$ above the graph.
The charts can be chosen so that the graphs over $B'_r$ cover
$\partial\Om_t$ for all sufficiently small $t$. In each chart,
\eqref{c2:eq:smoothing}, \eqref{c2:eq:nesting} and implicit
differentiation of $\rho_t(x',b_t(x'))=0$ give
\begin{equation}\label{c2:eq:identification}
 0\le b_t-b_{2t}\le Ct^3,\qquad
 |D'b_t-D'b_{2t}|\le Ct^2.
\end{equation}
Let $\nu_t$ be the interior unit normal to $\partial\Om_t$. Since
$u_t-\Phi_t$ vanishes on the boundary, its normal derivative satisfies
\[
 \begin{aligned}
 q_t(x')&:=\partial_{\nu_t}(u_t-\Phi_t)(x',b_t(x'))\\
        &=\sqrt{1+|D'b_t(x')|^2}\,\partial_n(u_t-\Phi_t)(x',b_t(x')) .
 \end{aligned}
\]

\begin{lemma}\label{c2:lem:comparison}
For all sufficiently small dyadic $t$,
\begin{align}
 \|u_t-u_{2t}\|_{L^\infty(\Om_t)}
     &\le C_\eps t^{2+\alpha-\eps},\label{c2:eq:successive}\\
 \|q_t-q_{2t}\|_{L^\infty(B'_{2r})}
     &\le C_\eps t^{1+\alpha/2-\eps}\label{c2:eq:normal-difference}
\end{align}
in each of the fixed charts.
\end{lemma}

\begin{proof}
On $\partial\Om_t$ the boundary displacement, the gradient bound
and \eqref{c2:eq:smoothing} give $|u_t-u_{2t}|\le Ct^3$.
Since $|d_t-d_{2t}|\le Ct^3$, \eqref{c2:eq:density} gives
\[
 |g_t-g_{2t}|\le Ct^\alpha\mathbf{1}_{\{d_t<C_1t\}}+Ct^4
 \qquad\text{in }\Om_t
\]
for a fixed $C_1$. For the linearization at either solution,
$\cL d_t\le-(\theta/2)\tr a\le-c$ in $\Om_t$.
On $d_t<2C_1t$, \eqref{c2:eq:ellipticity} and
\eqref{c2:eq:preliminary} give
$a^{ij}(d_t)_i(d_t)_j\ge c_\eps t^\eps\ge c'_\eps d_t^\eps$.
Thus \eqref{c2:eq:barrier} applies to both linearizations on
$\Om_t$. For $J=Ct^\alpha$ and $s=C_1t$, take
\[
 v=\psi_{J,s}(d_t)+C_2t^4d_t+C_3t^3 .
\]
For sufficiently large $C_2,C_3$, concavity gives
$\cL(u_{2t}-u_t)\ge g_{2t}-g_t\ge\cL v$
for the linearization at $u_t$. The maximum principle and the
corresponding inequality with $u_t,u_{2t}$ interchanged give
$|u_t-u_{2t}|\le v$. Now \eqref{c2:eq:barrier} proves
\eqref{c2:eq:successive}.

Fix $x'\in B'_{2r}$, put $P=(x',b_t(x'))$ and set
$h=t^{1+\alpha/2}$. For small $t$ the vertical segment from $P$ to
$P+he_n$ stays in the chart and lies in both closed domains. Taylor's formula and
\eqref{c2:eq:preliminary} give
\[
 |\partial_n(u_t-u_{2t})(P)|
 \le 2h^{-1}\|u_t-u_{2t}\|_{L^\infty(\Om_t)}+C_\eps ht^{-\eps}
 \le C_\eps t^{1+\alpha/2-\eps}.
\]
The segment from $(x',b_{2t}(x'))$ to $P$ also lies in
$\overline{\Om_{2t}}$. Along this segment $Du_{2t}$ varies by at most $C_\eps t^{3-\eps}$ by
\eqref{c2:eq:preliminary} and \eqref{c2:eq:identification}.
The corresponding first derivatives of $\Phi_t,\Phi_{2t}$ differ
by at most $Ct^2$, by \eqref{c2:eq:smoothing}. The factors
$\sqrt{1+|D'b_t|^2}$ and $\sqrt{1+|D'b_{2t}|^2}$ also differ by
at most $Ct^2$. The gradient bound gives \eqref{c2:eq:normal-difference}.
\end{proof}

\begin{lemma}\label{c2:lem:one-sided}
In each fixed chart, for every unit vector $e\in\R^{n-1}$,
\begin{equation}\label{c2:eq:one-sided}
 \partial_{ee}q_t\le C_\eps t^{-1-\eps}\qquad\text{in }B'_{2r}.
\end{equation}
\end{lemma}

\begin{proof}
First take coordinates centered at an arbitrary boundary point with
inward normal $e_n$, as in \eqref{c2:eq:rigid}. Write $x_n=b(x')$
for the graph in these coordinates, so that $b(0)=D'b(0)=0$.
Choose any unit tangent as $e_1$ and put $W=X^2w$, with $X,w$
as in \eqref{c2:eq:rigid}. Differentiating the equation twice
under the rigid motions generated by $X$ gives
$\cL(X^2u_t)\ge X^2g_t$ by concavity~\cite[Lemma~2.1]{ITW}. Thus the fourth derivative
bound on $\Phi_t$ and \eqref{c2:eq:density} give
\[
 \cL W\ge-J-Ct^{-1}\tr a\qquad(d_t<t),\qquad J=Ct^{\alpha-2}.
\]

On the graph write $\xi=e_1-bh$ and
$\delta=h\cdot x'-\xi\cdot D'b$. For every smooth $v$,
\[
 (Xv)(x',b(x'))=\xi\cdot D'\bigl(v(x',b(x'))\bigr)
                       +\delta v_n(x',b(x')) .
\]
In particular, $Xw=\delta w_n$ on the graph.
Since $h_i=b_{1i}(0)$, Taylor's formula yields
\[
 \delta=O(|x'|^2),\qquad
 \xi\cdot D'\delta=-\sum_{j<n}b_{11j}(0)x_j+
                  O\bigl((1+\|b\|_{C^4})|x'|^2\bigr).
\]
Applying the graph identity to $v=Xw$, we obtain
\[
 |W-\ell|\le C_\eps t^{-1}|x'|^2,\qquad
 \ell(x')=-w_n(0)\sum_{j<n}b_{11j}(0)x_j .
\]
Indeed $w_n$ varies by at most $C(1+M_t)|x'|$ along the
graph, the terms containing $\delta$ are bounded by
$C(1+M_t)|x'|^2$, and $\|b\|_{C^4}\le Ct^{-1}$.

Extend $\ell$ independently of $x_n$. In the boundary neighborhood used in the proof of
Lemma~\ref{c2:lem:preliminary} compare $W-\ell$ from above with
\[
 C_1t^{-1-\eps}d_t+C_2t^{-1}|x|^2+\psi_{J,t}(d_t).
\]
On the graph the preceding estimate gives $W-\ell\le C_\eps t^{-1}|x'|^2$.
On $\partial B_R$, $|W-\ell|\le C(1+M_t)\le C_\eps t^{-\eps}$ and
$R$ is fixed. Thus $C_2t^{-1}|x|^2$ bounds $W-\ell$ on both portions
of the boundary. On $d_t=t$ the first term suffices. After fixing $C_2$, take
$C_1$ large and use $\cL d_t\le-(\theta/2)\tr a$ and
\eqref{c2:eq:barrier} to obtain the differential inequality.
The maximum principle gives the comparison, with equality at zero.
Taking inward derivatives there gives
\[
 W_n(0)\le C_\eps(t^{-1-\eps}+Jt^{1-\eps})
          \le C_\eps t^{-1-\eps}.
\]
At zero $w_{i1}=-b_{i1}w_n$ for $i<n$, and differentiation of
\eqref{c2:eq:rigid} gives
\[
 W_n(0)=w_{11n}(0)+b_{11}(0)w_{nn}(0)+|h|^2w_n(0).
\]
This is the second $x_1$ derivative at zero of
$\sqrt{1+|D'b|^2}\,w_n(x',b(x'))$, the normal derivative in
the centered chart. The same upper bound holds in every unit
direction of its coordinate plane.

To return to a fixed chart at $x'$, project the curve
$s\mapsto(x'+se,b_t(x'+se))$ orthogonally onto the tangent plane
at $(x',b_t(x'))$. Its coordinates $y(s)$ in that plane satisfy
$y(0)=0$ and $|y'(0)|+|y''(0)|\le C$ by the uniform $C^2$
graph bounds. Write $q$ for the normal derivative in the centered
chart. Then $q_t(x'+se)=q(y(s))$, and $|D'q(0)|\le C_\eps
t^{-\eps}$ by \eqref{c2:eq:preliminary}. Applying the ordinary
chain rule twice gives \eqref{c2:eq:one-sided}.
\end{proof}

\begin{proposition}\label{c2:prop:boundary}
For sufficiently small $t_0$ and $t_j=2^{-j}t_0$,
\begin{equation}\label{c2:eq:boundary-bound}
 \sup_j\|D^2u_{t_j}\|_{L^\infty(\partial\Om_{t_j})}\le C .
\end{equation}
\end{proposition}

\begin{proof}
Take $t_0<\min\{r,1\}/2$ and set $h=t^{1+\alpha/4}$.
For $x'\in B'_{r+t}$ and $|e|=1$, the segment from $x'-he$ to
$x'$ lies in $B'_{r+2t}$, since $h\le t$. By
\eqref{c2:eq:one-sided} and Taylor's formula,
\[
 \partial_e q_t(x')
 \le \frac{q_t(x')-q_t(x'-he)}{h}+C_\eps t^{-1-\eps}h .
\]
Replacing $q_t$ by $q_{2t}$ changes the difference quotient by at most $2C_\eps t^{1+\alpha/2-\eps}/h$ by
\eqref{c2:eq:normal-difference}. The resulting quotient is at
most $\|D'q_{2t}\|_{L^\infty(B'_{r+2t})}$. Taking the supremum
over $e$ and $x'$ gives
\begin{equation}\label{c2:eq:mixed-sum}
 \|D'q_t\|_{L^\infty(B'_{r+t})}
 \le \|D'q_{2t}\|_{L^\infty(B'_{r+2t})}
     +C_\eps t^{\alpha/4-\eps}.
\end{equation}
Here $\alpha/4-\eps>0$, so the errors are summable over dyadic $t$.
The initial norm is bounded by \eqref{c2:eq:preliminary}. We obtain $|D'q_t|\le C$
on $B'_r$ in each chart. Differentiating the normal derivative on
the graph gives
\[
 \partial_i q_t
 =\bigl(e_i+(b_t)_i e_n\bigr)\cdot D^2(u_t-\Phi_t)\nu_t,
 \qquad i<n,
\]
since $D(u_t-\Phi_t)$ is normal and the derivative of the unit
normal is perpendicular to it. The smaller charts cover
$\partial\Om_t$, and their slopes are uniformly bounded.
This bounds all mixed derivatives there. Proposition~\ref{bdy:margin}
and \eqref{c2:eq:block} give \eqref{c2:eq:boundary-bound}.
\end{proof}

\subsection{Global second derivative estimates}\label{c2:global-trace}

\begin{proposition}\label{c2:prop:trace}
For a fixed $C_0$ and every sufficiently small dyadic $t$,
\begin{equation}\label{c2:eq:trace-sum}
 M_t\le\max\{C_0,M_{2t}\}+Ct^{\alpha/10}.
\end{equation}
Consequently $\sup_j M_{2^{-j}t_0}\le C$.
\end{proposition}

\begin{proof}
Choose $C_0\ge2$ so that $S_t\le C_0$ on $\partial\Om_t$ by \eqref{c2:eq:boundary-bound},
and put $e_t=C_\eps t^{\alpha-\eps}$ as in \eqref{c2:eq:localization}.
It suffices to consider $M_t>C_0+e_t$. Then there is a point $x_t$ with
$d_t(x_t)\ge t$ and $S_t(x_t)\ge M_t-e_t$. Choose $R=ct$ so that
$B_R(x_t)\Subset\Om_t$, and put $h=t^{1+\alpha/4}$. By
\eqref{c2:eq:preliminary}, $M_t^{k/2}h/R\le C_\eps t^{\alpha/4-k\eps/2}$
and $e_t\le M_t/8$ for all sufficiently small $t$.

Let $K$ be the nonnegative kernel in Lemma~\ref{c2:lem:persistence}
on $B_R(x_t)$ at scale $h$. Since $\int K\le1$ and $\Om_t\subset\Om_{2t}$, the
estimate in that lemma applied to $w=u_t-u_{2t}$ gives
\[
 \int K\Delta u_t
 =\int K\Delta u_{2t}+\int K\Delta(u_t-u_{2t})
 \le M_{2t}+Ch^{-2}\|u_t-u_{2t}\|_{L^\infty(\Om_t)} .
\]
The lower bound \eqref{c2:eq:trace-average}, together with
\eqref{c2:eq:preliminary} and \eqref{c2:eq:successive}, now gives
\[
 M_t\le M_{2t}
  +C_\eps\left(t^{\alpha-\eps}+t^{\alpha/2-\eps}
                         +t^{\alpha/4-(1+k/2)\eps}\right).
\]
Indeed the comparison term is $O(t^{\alpha/2-\eps})$, and the other
losses are $e_t$, $CM_th^{\alpha/2}$ and $CM_t^{1+k/2}h/R$.
The smallest exponent is
$\alpha/4-(1+k/2)\eps=49\alpha/200>\alpha/10$. If $M_t\le C_0+e_t$, then
\eqref{c2:eq:trace-sum} follows directly from the bound on $e_t$. Summing on dyadic scales and using
\eqref{c2:eq:preliminary} at the initial scales bounds $M_t$ uniformly.
\end{proof}

\subsection{Schauder estimates}\label{c2:holder}

\begin{proof}[Proof of Theorem~\ref{thm:c2alpha}]
Proposition~\ref{c2:prop:trace} gives $|D^2u_t|\le\Delta u_t\le C$
at all dyadic scales. Since $G(D^2u_t)\ge m$, the Hessians lie in
the interior of the compact convex set
\[
 \cK=\{H\in\Ga:G(H)\ge m/2,\ \tr H\le 2C\}\Subset\Ga .
\]
The eigenvalues of $(G^{ij})$ have uniform positive lower and finite upper
bounds on $\cK$. Thus $\widehat G(B):=\inf_{H\in\cK}G^{ij}(H)B_{ij}$
is globally uniformly elliptic and concave. By homogeneity and concavity,
$\widehat G=G$ on $\cK$ and $\widehat G(0)=0$.
The global estimate of Silvestre and
Sirakov~\cite[Theorem~1.4(b)]{SS}, applied to the convex operator
$B\mapsto-\widehat G(-B)$ and its solution $-u_t$, with right-hand
side $-g_t$ and boundary values $-\Phi_t$, gives
$\|u_t\|_{C^{2,\gamma}(\overline{\Om_t})}\le C$
for some $0<\gamma<\alpha$. Here the $C^{2,\alpha}$ norms of the
boundary charts and boundary values are uniform.

To recover the exponent $\alpha$, fix $x_0$ in an interior or boundary
patch of radius $r$ and put $a_0^{ij}=G^{ij}(D^2u_t(x_0))$. For any two Hessians
in the patch, $|G^{ij}-a_0^{ij}|\le Cr^\gamma$ on the
segment joining them. This follows from the $C^{2,\gamma}$ estimate and
the smoothness of $G$ on $\cK$. Applying the mean value formula to
$G(D^2u_t)=g_t$ gives, with seminorms taken on the patch,
\[
 [a_0^{ij}(u_t)_{ij}]_{C^\alpha}
 \le [g_t]_{C^\alpha}+Cr^\gamma[D^2u_t]_{C^\alpha}.
\]
Apply the interior Schauder estimate for $a_0^{ij}\partial_{ij}$ on
the patch of half the radius. On a boundary patch apply the boundary
estimate to $u_t-\Phi_t$, which has zero boundary values. The term
$a_0^{ij}(\Phi_t)_{ij}$ has uniformly bounded $C^\alpha$ norm.
In the rescaled estimates \cite[Chapter~6]{GT}, the two seminorms
have the same factor $r^{2+\alpha}$ and the boundary charts have
uniformly bounded $C^{2,\alpha}$ norms. Thus the coefficient of the last term is at most
$C_1r^\gamma$, with $C_1$ independent of $r,t$. The other terms are
bounded for fixed $r$ by the $C^2$ estimate and the norms of the data.
Any two points at distance less than $cr$ lie in one of the smaller
patches, for a uniform $c>0$. For the other pairs the $C^2$ bound suffices.
Taking the global seminorm and choosing $r$ so that $C_1r^\gamma\le1/2$
gives the uniform $C^{2,\alpha}$ estimate by absorption.

By \eqref{c2:eq:nesting} the restrictions to $\overline\Om$ are
uniformly bounded in $C^{2,\alpha}(\overline\Om)$. A subsequence
therefore converges in $C^2(\overline\Om)$.
The limit inherits the Hessian H\"older bound and takes the boundary
values $\varphi$,
by the boundary displacement, the gradient bound and $\Phi_t\to\Phi$.
Its Hessian lies in $\cK$ and its density is $f$ by
\eqref{c2:eq:density}. Comparison for Hessian measures
\cite[Theorem~3.1]{TW1} identifies it with the prescribed solution.
\end{proof}

\section*{Acknowledgements}

The author is deeply grateful to Professors Dongsheng Li and Yu Yuan for their guidance and support.
AI tools were used as an aid in preliminary mathematical exploration and in refining the exposition.
All mathematical arguments were independently verified by the author, who takes full responsibility for the manuscript.

\begingroup
\renewcommand{\baselinestretch}{1}\selectfont

\endgroup


\begin{thebibliography}{99}

\bibitem{Caffarelli}
L. A. Caffarelli,
\emph{Interior $W^{2,p}$ estimates for solutions of the Monge--Amp\`ere equation},
Ann. of Math. (2) \textbf{131} (1990), no.~1, 135--150.

\bibitem{CNS}
L. Caffarelli, L. Nirenberg and J. Spruck,
\emph{The Dirichlet problem for nonlinear second-order elliptic equations.
III. Functions of the eigenvalues of the Hessian},
Acta Math. \textbf{155} (1985), no.~3--4, 261--301.

\bibitem{CWW}
C. Chen, X.-J. Wang and Y. Wu,
\emph{$C^{1,\alpha}$ regularity of convex hypersurfaces with prescribed
curvature measures},
Indiana Univ. Math. J. \textbf{72} (2023), no.~1, 331--351.

\bibitem{CZZ}
R. Chen, X. Zhou and R. Zhu,
\emph{Interior $C^{2,\alpha}$ regularity for the quadratic Hessian equation},
preprint, \href{https://arxiv.org/abs/2608.29484v2}{arXiv:2608.29484v2}, 2026.

\bibitem{CW}
K.-S. Chou and X.-J. Wang,
\emph{A variational theory of the Hessian equation},
Comm. Pure Appl. Math. \textbf{54} (2001), no.~9, 1029--1064.

\bibitem{DK}
S. Dinew and S. Ko\l odziej,
\emph{A priori estimates for complex Hessian equations},
Anal. PDE \textbf{7} (2014), no.~1, 227--244.

\bibitem{Garding}
L. G{\aa}rding,
\emph{An inequality for hyperbolic polynomials},
J. Math. Mech. \textbf{8} (1959), no.~6, 957--965.

\bibitem{GT}
D. Gilbarg and N. S. Trudinger,
\emph{Elliptic Partial Differential Equations of Second Order},
reprint of the 1998 edition, Classics in Mathematics,
Springer, Berlin, 2001.

\bibitem{Guan}
B. Guan,
\emph{The Dirichlet problem for a class of fully nonlinear elliptic equations},
Comm. Partial Differential Equations \textbf{19} (1994), no.~3--4, 399--416.

\bibitem{He}
R. He,
\emph{Gradient estimates for Hessian equations with measure data},
preprint, \href{https://arxiv.org/abs/2610.38846v1}{arXiv:2610.38846v1}, 2026.

\bibitem{Ivochkina}
N. M. Ivochkina,
\emph{Solution of the Dirichlet problem for some equations of
Monge--Amp\`ere type},
Math. USSR-Sb. \textbf{56} (1987), no.~2, 403--415.

\bibitem{ITW}
N. M. Ivochkina, N. S. Trudinger and X.-J. Wang,
\emph{The Dirichlet problem for degenerate Hessian equations},
Comm. Partial Differential Equations \textbf{29} (2004), no.~1--2,
219--235.

\bibitem{KrylovBoundary}
N. V. Krylov,
\emph{Boundedly nonhomogeneous elliptic and parabolic equations in a domain},
Math. USSR-Izv. \textbf{22} (1984), no.~1, 67--97.

\bibitem{LMZ}
H. Liu, S. Ma and S. Zang,
\emph{Gradient potential estimates for $k$-Hessian equations},
preprint, \href{https://arxiv.org/abs/2611.08184v1}{arXiv:2611.08184v1}, 2026.

\bibitem{LyuInterior}
Y. Lyu,
\emph{Interior $W^{2,p}$ regularity for the sigma-$2$ equation},
preprint, \href{https://arxiv.org/abs/2609.08297v1}{arXiv:2609.08297v1}, 2026.

\bibitem{LyuK2}
Y. Lyu,
\emph{Global $W^{2,p}$ and $C^{2,\alpha}$ regularity for the sigma-$2$ equation},
preprint, \href{https://arxiv.org/abs/2609.11909v1}{arXiv:2609.11909v1}, 2026.

\bibitem{LyuEps}
Y. Lyu,
\emph{$C^{1,\alpha}$ and $W^{2,1+\varepsilon}$ regularity for the
$2$-Hessian equation},
preprint, \href{https://arxiv.org/abs/2610.35839v1}{arXiv:2610.35839v1}, 2026.

\bibitem{LyuWone}
Y. Lyu,
\emph{Interior $W^{2,1}$ regularity for the $k$-Hessian equation},
manuscript.

\bibitem{MooneySingular}
C. Mooney,
\emph{$W^{2,1}$ estimate for singular solutions to the Monge--Amp\`ere equation},
Ann. Sc. Norm. Super. Pisa Cl. Sci. (5) \textbf{14} (2015), no.~4, 1283--1303.

\bibitem{Reilly}
R. C. Reilly,
\emph{On the Hessian of a function and the curvatures of its graph},
Michigan Math. J. \textbf{20} (1973), no.~4, 373--383.

\bibitem{Savin}
O. Savin,
\emph{Global $W^{2,p}$ estimates for the Monge--Amp\`ere equation},
Proc. Amer. Math. Soc. \textbf{141} (2013), no.~10, 3573--3578.

\bibitem{SavinPoint}
O. Savin,
\emph{Pointwise $C^{2,\alpha}$ estimates at the boundary for the
Monge--Amp\`ere equation},
J. Amer. Math. Soc. \textbf{26} (2013), no.~1, 63--99.

\bibitem{SS}
L. Silvestre and B. Sirakov,
\emph{Boundary regularity for viscosity solutions of fully nonlinear elliptic equations},
Comm. Partial Differential Equations \textbf{39} (2014), no.~9, 1694--1717.

\bibitem{Trudinger1990}
N. S. Trudinger,
\emph{The Dirichlet problem for the prescribed curvature equations},
Arch. Rational Mech. Anal. \textbf{111} (1990), no.~2, 153--179.

\bibitem{Trudinger1995}
N. S. Trudinger,
\emph{On the Dirichlet problem for Hessian equations},
Acta Math. \textbf{175} (1995), no.~2, 151--164.

\bibitem{TW1}
N. S. Trudinger and X.-J. Wang,
\emph{Hessian measures I},
Topol. Methods Nonlinear Anal. \textbf{10} (1997), no.~2, 225--239.

\bibitem{TW2}
N. S. Trudinger and X.-J. Wang,
\emph{Hessian measures II},
Ann. of Math. (2) \textbf{150} (1999), no.~2, 579--604.

\bibitem{TWMA}
N. S. Trudinger and X.-J. Wang,
\emph{Boundary regularity for the Monge--Amp\`ere and affine maximal
surface equations},
Ann. of Math. (2) \textbf{167} (2008), no.~3, 993--1028.

\bibitem{Urbas}
J. Urbas,
\emph{On the existence of nonclassical solutions for two classes of fully
nonlinear elliptic equations},
Indiana Univ. Math. J. \textbf{39} (1990), no.~2, 355--382.

\bibitem{Urbas2000}
J. Urbas,
\emph{Some interior regularity results for solutions of Hessian equations},
Calc. Var. Partial Differential Equations \textbf{11} (2000), no.~1, 1--31.

\bibitem{Urbas2001}
J. Urbas,
\emph{An interior second derivative bound for solutions of Hessian equations},
Calc. Var. Partial Differential Equations \textbf{12} (2001), no.~4, 417--431.

\bibitem{Wang}
X.-J. Wang,
\emph{Regularity for Monge--Amp\`ere equation near the boundary},
Analysis \textbf{16} (1996), no.~1,
\href{https://doi.org/10.1524/anly.1996.16.1.101}{101--107}.
\bibitem{ZZ}
X. Zhou and R. Zhu,
\emph{Interior $C^{2,\alpha}$ regularity for convex solutions of the
2-Hessian equation},
preprint, \href{https://arxiv.org/abs/2608.24604v1}{arXiv:2608.24604v1}, 2026.
\end{thebibliography}
\end{document}